\documentclass[3p,times,chicago]{elsarticle}

\usepackage{amssymb}
\oddsidemargin \dimexpr 1in -1in
\topmargin \dimexpr 1in -1in
\textwidth \dimexpr \pdfpagewidth -2\oddsidemargin -2in
\textheight \dimexpr \pdfpageheight -2\topmargin -2in

\usepackage{amsmath,amssymb,amsthm,textcomp}
\usepackage{amsfonts,graphicx}
\usepackage[mathscr]{eucal}
\usepackage{color}
\usepackage{hyperref}
\usepackage[table]{xcolor}
\usepackage{booktabs}

\usepackage[]{natbib}
\usepackage{subfig}
\usepackage{epsfig}
\usepackage{caption}
\theoremstyle{definition}
\usepackage{float}
\usepackage{mathtools}  
\usepackage{diffcoeff}  
\floatstyle{plaintop}
\usepackage[figuresright]{rotating}
\usepackage{natbib}
\usepackage{amssymb,amsmath, amsfonts}
\usepackage{enumerate}
\usepackage{amsthm}
\usepackage[utf8]{inputenc}
\usepackage[english]{babel}
\setcitestyle{authoryear,open={(},close={)}} 

\usepackage{graphicx}
\usepackage{caption}
\newtheorem{theorem}{Theorem}[section]
\newtheorem{proposition}{Proposition}[section]

\newtheorem{remark}{Remark}[section]
\newtheorem{lemma}{Lemma}[section]

\numberwithin{equation}{section}

\newtheoremstyle
{remarkstyle}
{}
{11pt}
{}
{}
{\bfseries}
{:}
{     }
{\thmname{#1} \thmnumber{#2} }

\theoremstyle{remarkstyle}

\makeatletter
\def\ps@pprintTitle{%
	\let\@oddhead\@empty
	\let\@evenhead\@empty
	\let\@oddfoot\@empty
	\let\@evenfoot\@oddfoot
}
\makeatother
\begin{document}
	
	\begin{frontmatter}
		
		
		

\title{Generalized Fractional Risk Process}

\author[label1]{Ritik Soni}

\author[label1]{Ashok Kumar Pathak \corref{cor1}}
\ead{ashokiitb09@gmail.com}
\address[label1]{Department of Mathematics and Statistics, Central University of Punjab, Bathinda, India}
\cortext[cor1]{Corresponding author}

\begin{abstract} In this paper, we define a compound generalized fractional counting process (CGFCP) which is a generalization of the compound versions  of  several well-known fractional counting processes.  We obtain its mean, variance, and the fractional differential equation governing the probability law. Motivated by \cite{Kumar2020}, we introduce a fractional risk process by considering CGFCP as the surplus  process and call it generalized fractional risk process (GFRP). We study the martingale property of the GFRP and show that GFRP and the associated increment process exhibit the long-range dependence (LRD) and the short-range dependence (SRD) property, respectively. We also define an alternative to GFRP, namely AGFRP  which is premium wise different from the GFRP. Finally, the asymptotic structure of the ruin probability for the AGFRP is established in case of light-tailed and heavy-tailed claim sizes.
\end{abstract}

\end{frontmatter}


	
\section{Introduction}
In recent years, researchers have shown a considerable interest in studying the time changed variants of the Poisson process.   Fractional generalizations of the Poisson processes are broadly classified into three categories: time-fractional, space-fractional, and space-time fractional (see \cite{Beghin2009}, \cite{Laskin2009}, \cite{Meerschaert2011}, and \cite{Orsingher2012}). Potential applications of these processes are  observed in the areas of finance, actuarial science, economics, biology, hydrology, and reliability (see \cite{Doukhan2002}, \cite{Biard2014}, \cite{Kumar2020}, \cite{Crescenzo2023}, and \cite{Ritik2024}). The compound version of the Poisson process generalizes the various unit jump processes to a random jump.  It can accommodate time dependencies and complex event patterns that arise from a variety of real-world phenomena. \cite{meerschaert2004limit} first considered the fractional generalization of the compound Poisson process (CPP). Later on, alternative forms of the fractional compound Poisson processes have been explored by \cite{Beghin2009, beghin2012alternative, Scalas2012}. \cite{beghin2014fractional} proposed two fractional versions of nonnegative, integer-valued compound Poisson processes and derived their governing fractional  differential equations. Further, \cite{Biard2014} and \cite{Maheshwari2016} examined the long-range dependence (LRD) property of the fractional compound Poisson processes and showed their applicability in ruin theory. Other notable references in this direction are \cite{gupta2023fractional}, \cite{Kumar2020},  \cite{sengar2020subordinated}, and \cite{Kataria2022}.

The generalized fractional counting process (GFCP) is a unified extension of the fractional Poisson processes (FPPs). It can accommodate  random jumps with varying intensity and could be useful in studying and explaining  various real-world complex phenomena. \cite{Crescenzo2015} presented a subordinated form of the GFCP and established its connection with the fractional Cauchy problem. Further,  \cite{Kataria2021} discussed its distributional properties and dependence  structure. Recently, \cite{khandakar2023time} introduced a mixed fractional counting process and applied it to risk theory.

In this paper, we present a compound generalized fractional counting process that extends the compound versions of several well-known fractional counting processes (for example, see \cite{Scalas2012}, \cite{Kataria2021}, and \cite{Kataria2022}). We consider a generalized fractional risk process (GFRP) based on the GFCP and discuss its martingale characterization. The GFRP  generalizes the fractional risk models studied by \cite{Kumar2020}. It can also be applied to various classical risk models based on the GCP and convoluted Poisson process which have not been studied in the literature to the best of  our  knowledge.  We have examined the dependence properties of the GFRP and its associated increment process.  Additionally, we define AGFRP, a premium-wise alternative to GFRP. Finally, for both light-tailed and heavy-tailed claim sizes, the asymptotic structure of the ruin probability for the AGFRP is established.

The paper is structured as follows: In Section 2, we present some preliminary notations and definitions. In Section 3, we define the CGFCP and discuss its distributional properties. We consider the GFRP and examine its main characteristics in Section 4. Finally, an alternative to the GFRP is introduced and its asymptotic structure is explored in Section 5.

\section{Preliminaries}
\noindent Here, we present some known results which are used in the subsequent sections.
Let $\mathbb{N}$ be the set of natural numbers and $\mathbb{N}_0 = \mathbb{N} \cup \{0\}$. Let $\mathbb{R}$ and $\mathbb{C}$ denote the set of real numbers and set of complex numbers, respectively. Let $x \wedge y = \min\{x,y\}$ and $x \vee y = \max\{x,y\}.$
\subsection{Special Functions and Generalized Derivatives}\hfill \\
(i) The generalized Mittag-Leffler function $E_{\beta, \gamma}^{\alpha}(z)$ is defined as (see \cite{Podlubny})
\begin{equation*}
E_{\beta, \gamma}^{\alpha}(z) = \sum_{k=0}^{\infty} \frac{z^k}{k!\Gamma(\gamma+\beta k)}\frac{\Gamma(\alpha +k)}{\Gamma(\alpha)}, \;\;\beta, \gamma, \alpha,  z \in \mathbb{C} \text{ and } \text{Re}(\beta)>0, \text{Re}(\gamma)>0, \text{Re}(\alpha)>0.
\end{equation*}
(ii) Let $f :[a,b] \subset \mathbb{R}\longrightarrow\mathbb{R}$ be $(n+1)$-times continuous differentiable function for $n < \tau <n+1$. Then, the Riemann-Liouville (RL) fractional derivative of order $\tau>0$ is defined as (see \cite{Podlubny})
\begin{equation*}
_aD^\tau_tf(t)=\bigg( \frac{d}{dt}\bigg)^{n+1}\int_{a}^{t}(t-u)^{n-\tau}f(u)du.
\end{equation*}
\subsection{L\'{e}vy Subordinator}\hfill \\
A L\'{e}vy subordinator denoted by $\{S_\psi(t)\}_{t \geq 0}$ is a  non-decreasing L\'{e}vy process with Laplace transform (see \cite[Section 1.3.2]{Applebaum})
\begin{equation*}
\mathbb{E}\left(e^{-uS_\psi(t)}\right) = e^{-t \psi(u)},\;\; u \geq 0,
\end{equation*}
where $\psi(u)$ is Laplace exponent given by (see \cite[Theorem 3.2]{Schilling})
\begin{equation*}
\psi(u) = \eta u +\int_{0}^{\infty} (1-e^{-ux}) \nu (dx), \;\; \eta \geq 0.
\end{equation*}
Here $\eta$ is the drift coefficient and $\nu$ is a non-negative L\'{e}vy measure on the positive half-line satisfying
\begin{equation*}
\int_{0}^{\infty} \min \{x,1\} \nu (dx) < \infty, \;\;\; \text{ and } \;\;\; \nu ([0,\infty)) =\infty,
\end{equation*}
so that $\{S_\psi(t)\}_{t \geq 0}$ has strictly increasing sample paths almost surely (a.s.) (for more details see  \cite[Theorem 21.3]{Sato}).
\subsection{Inverse Subordinator}
Let $\{S_\psi(t)\}_{t \geq 0}$ be a stochastically continuous L\'{e}vy subordinator with independent and stationary increments. Then, the inverse subordinator denoted by $\{\mathcal{W}(t)\}_{t \geq 0}$ is defined as the first passage time of $S_\psi$, i.e.
\begin{equation*}
\mathcal{W}(t) = \inf\{r \geq0 : S_\psi(r) > t \}.
\end{equation*}
The process $\{\mathcal{W}(t)\}_{t \geq 0}$ is non-markovian with non-stationary and non-independent increments which satisfies
\begin{equation*}
\text{Pr}\{\mathcal{W}(t) > r\} \leq \text{Pr}\{S_\psi(t) \leq t\}.
\end{equation*}
\subsubsection{Inverse $\beta$-Stable Subordinator}
For $\beta \in (0,1),$ let $\{S_\beta (t)\}_{t \geq 0}$ be $\beta$-stable subordinator with Laplace transform given by $\mathbb{E}[e^{-uS_\beta (t)}] = e^{-tu^{\beta}}, \; u > 0$. Then, the inverse $\beta$-stable subordinator $\{\mathcal{W}_\beta(t)\}_{t \geq 0}$  is defined by
\begin{equation*}
\mathcal{W}_\beta(t) = \inf\{r \geq0 : S_\beta(r) > t \}.
\end{equation*}
The mean of the $\mathcal{W}_\beta(t)$ is given by (see \cite{Leonenko2014})
\begin{equation*}
\mathbb{E}[\mathcal{W}_\beta(t)] = \frac{t^{\beta}}{\Gamma(\beta +1)}.
\end{equation*}
Let $B(\beta, \beta+1)$ and $B(\beta, \beta+1; s/t)$ be the beta function and the incomplete beta function, respectively. It is known that for $0 <s \leq t$, the covariance of $\mathcal{W}_\beta(t)$ is expressed as (see \cite[formula (10)]{Leonenko2014})
\begin{align}\label{asy11}
\text{Cov}[\mathcal{W}_\beta(s), \mathcal{W}_\beta(t)] 
&= \frac{1}{\Gamma^2(\beta+1)} \left(\beta s^{2 \beta} B(\beta, \beta+1) +\beta t^{2\beta}B(\beta, \beta+1; s/t)-(ts)^{\beta}\right),\nonumber\\
& \sim \frac{1}{\Gamma^2(\beta+1)} \left(\beta s^{2 \beta} B(\beta, \beta+1) - \frac{\beta^2 s^{\beta+1}}{(\beta+1)t^{1-\beta}}\right),
\end{align}
for a fixed $s$ and large $t$.
\subsection{Dependence Structure} 
For $0 <s <t,$ a stochastic process $\{X(t)\}_{t \geq 0}$ is said to have LRD property if 
\begin{equation*}
\text{ Corr}[X(s), X(t)]\sim c(s){t^{-d}}\; \text{as}~  t\rightarrow \infty,
\end{equation*}
where the constant $c(s)$ is depending on $s$ and $d \in (0,1)$. For $d \in (1,2)$, the process exhibits SRD property.

\subsection{Generalized Fractional Counting Process}
For  $0 < \beta \leq 1$, let $\{\mathcal{N}^{\beta}(t)\}_{t \geq 0}$ be a generalized fractional counting process (GFCP) which performs $k$ kinds of jump of amplitude $1,2,\ldots, k$ with positive rates $\lambda_1, \lambda_2,\ldots, \lambda_k,$ respectively for a fixed integer $k  \geq 1$. Its pmf $p_n^{\beta}(t) =$ Pr$\{\mathcal{N}^{\beta}(t) =n\} $ is given by (see \cite{Crescenzo2015} and \cite[formula (2)]{  Kataria2022})
\begin{equation*}
p_n^{\beta}(t) = \sum_{r=0}^{n}\sum_{\substack{
		j_1+ j_2+ \cdots+ j_n =r \\ \\ j_1+ 2j_2+ \cdots+ kj_k =n}} \frac{r! \lambda_1^{j_1}\lambda_2^{j_2}\cdots \lambda_k^{j_k}}{j_1!j_2! \cdots j_k!} t^{\beta r} E_{\beta, \beta r+1}^{r+1}(-\Delta t^{\beta}), \;\; n \geq 0,
\end{equation*}
where $\Delta = \lambda_1+\lambda_2+\cdots+\lambda_k$ and $j_1, j_2, \ldots, j_k$ are non-negative integers.\\
The governing fractional differential equation of order $\beta$ is 
\begin{equation}\label{de1}
\frac{d^{\beta}}{dt^{\beta}} 	p_n^{\beta}(t) = -\Delta 	p_n^{\beta}(t) + \sum_{j=1}^{n \wedge k}\lambda_j 	p_{n-j}^{\beta}(t), \;\; n \geq 0,
\end{equation}
where  $\frac{d^{\beta}}{dt^{\beta}}$ denotes the RL fractional differential operator.\\
The mean and variance of the GFCP respectively are given by 
\begin{equation}\label{mv1}
\mathbb{E}[\mathcal{N}^{\beta}(t)] = St^{\beta}, \;\;\; \text{Var}[\mathcal{N}^{\beta}(t)] = Rt^{2 \beta} + Tt^{\beta},
\end{equation}
where $S = \frac{\sum_{j=1}^k j \lambda_j}{\Gamma(\beta +1)}, \; R = \left(\frac{2}{\Gamma(2 \beta +1)}- \frac{1}{\Gamma^2(\beta+1)}\right)\left(\sum_{j=1}^{k}j \lambda_j\right)^2$ and $ T= \frac{\sum_{j=1}^k j^2 \lambda_j}{\Gamma(\beta +1)}$.

Let $\mathcal{N}(t)$ be a generalized counting process. Then, the GFCP can be viewed as a time-changed variant of the generalized counting process (GCP) of the form $$\mathcal{N}^{\beta}(t) := \mathcal{N}(\mathcal{W}_\beta(t)). $$
It may be seen that, when $\beta=1$, GFCP reduces to GCP (see \cite{Kataria2022} and \cite{Kataria2022a}).

\noindent For $0 < s \leq t$, the covariance function of the GFCP is given as  (see \cite[formula (24)]{Kataria2022})
\begin{align}\label{cov11}
\text{Cov}[\mathcal{N}^{\beta}(s), \mathcal{N}^{\beta}(t)] 
&= Ts^{\beta} + \left(\sum_{j=1}^{k}j \lambda_j\right)^2 \text{Cov}[\mathcal{W}_\beta(s), \mathcal{W}_\beta(t)]\\
& \sim Ts^{\beta} + S^2\left(\beta s^{2 \beta} B(\beta, \beta+1) - \frac{\beta^2 s^{\beta+1}}{(\beta+1)t^{1-\beta}}\right), \; \text{as}\; t \rightarrow \infty. \nonumber
\end{align}
\section{Compound Generalized Fractional Counting Process}
\noindent 	For $0 < \beta \leq 1$, let $\{\mathcal{N}^{\beta}(t)\}_{t \geq 0}$  be a GFCP with positive rates $\lambda_1, \lambda_2,\ldots, \lambda_k.$ We define the compound generalized fractional counting process (CGFCP) of the form
\begin{equation*}
C^{\beta}(t) = \sum_{j=1}^{\mathcal{N}^{\beta}(t)} X_j,\;\; t \geq 0,
\end{equation*}
where $X_j, \; j=1,2,\ldots$ are independent and identically distributed (i.i.d.) random variables with common cumulative distribution function (CDF) $H_X$ and are independent of $\mathcal{N}^{\beta}$. 

For $k=1$, the CGFCP reduces to the compound fractional Poisson process as GFCP reduces to fractional Poisson process (see \cite{Scalas2012}). By choosing $\lambda_j = \lambda$, the CGFCP can be viewed as compound representation of the fractional Poisson process of order $k$. For $\lambda_j = \lambda (1-\rho)\rho^{j-1}/(1-\rho^k), \; 0 \leq \rho < 1, $ the CGFCP leads to the compound P\'{o}lya-Aeppli process of order $k$ (see \cite{Kataria2022}). It is also noted that for a limiting case $k \rightarrow \infty $ with some suitable assumptions, the CGFCP is same as the compound convoluted fractional Poisson process (for more details see \cite{Kataria2021}).

Next, we present the state probabilities of the CGFCP.

\noindent Let $p_j =\; $Pr$\{X_1 =j\}\;$ for all $j \geq 1$ and  $p_n(i) = \;$Pr$\{X_1+X_2+\cdots+X_i = n\}$. Then, the pmf  of the CGFCP is given by

\begin{equation}\label{sp1}
q_n^\beta(t) = \sum_{i=1}^{n}r_n(i)p_i^{\beta}(t), \; \text{ with } \;\; q_0^\beta(t) = p_0^\beta(t).
\end{equation}
Now, we derive the governing fractional differential equation.

\begin{proposition}
	The state probabilities in (\ref{sp1}) satisfies the following fractional differential equation
	\begin{equation*}
	\diffp{^\beta}{t^{\beta}} q_n^\beta(t) = -\Delta q_n^\beta(t) + \sum_{j=1}^{ i\wedge k}\lambda_j  \sum_{l=1}^{n} r_l q_{n-l}^{\beta}(t). 
	\end{equation*}
\end{proposition}
\begin{proof}
	With the help of (\ref{de1}) and (\ref{sp1}), we have
	\begin{align}\label{stateprobability}
	\diffp{^\beta}{t^{\beta}} q_n^\beta(t) 
	&=  \sum_{i=1}^{n}r_n(i) \diffp{^\beta}{t^{\beta}}p_i^{\beta}(t)
	=  \sum_{i=1}^{n}r_n(i) \left(-\Delta 	p_i^{\beta}(t) + \sum_{j=1}^{ i\wedge k}\lambda_j 	p_{i-j}^{\beta}(t)\right)\nonumber\\
	&= -\Delta q_n^\beta(t) + \sum_{i=1}^{n}r_n(i) \sum_{j=1}^{ i\wedge k}\lambda_j 	p_{i-j}^{\beta}(t)
	= -\Delta q_n^\beta(t) + \sum_{j=1}^{ i\wedge k}\lambda_j \left[\sum_{i=1}^{n}r_n(i) p_{i-j}^{\beta}(t)\right]\nonumber\\
	&= -\Delta q_n^\beta(t) + \sum_{j=1}^{ i\wedge k}\lambda_j \left[\sum_{i=1}^{n}\left(\sum_{l=1}^n r_{n-l}(i-j)r_l\right)p_{i-j}^{\beta}(t)\right]\nonumber\\
	&= -\Delta q_n^\beta(t) + \sum_{j=1}^{ i\wedge k}\lambda_j \left[\sum_{l=1}^n r_l \sum_{i=1}^{n} r_{n-l}(i-j)p_{i-j}^{\beta}(t)\right].
	\end{align}
	Note that $r_k(0) = r_0(k)= 0$ for all $k \geq 1$, $r_n(k) =0$ for $n < k$, and $r_0(0)=1$. Therefore, we get
	\begin{align*}
	\sum_{l=1}^n r_l \sum_{i=1}^{n} r_{n-l}(i-j)p_{i-j}^{\beta}(t) 
	&= \sum_{l=1}^{n-1} r_l \sum_{i=1}^{n} r_{n-l}(i-j)p_{i-j}^{\beta}(t) +  r_n \sum_{i=1}^{n} r_{0}(i-j)p_{i-j}^{\beta}(t)\\
	&= \sum_{l=1}^{n-1} r_l \sum_{i=1}^{n-l} r_{n-l}(i)p_{i}^{\beta}(t) + r_n p_0^{\beta}(t)\\
	&= \sum_{l=1}^{n-1} r_l q_{n-l}^{\beta}(t) + r_n q_0^{\beta}(t)
	= \sum_{l=1}^{n} r_l q_{n-l}^{\beta}(t).
	\end{align*}
	Using this in (\ref{stateprobability}),  we get the required differential equation.
\end{proof}
 The mean and variance of the CGFCP are derived as
\begin{equation*}
\mathbb{E}[C^{\beta}(t)] = \mathbb{E}[\mathcal{N}^{\beta}(t)]\mathbb{E}[X_1] = St^{\beta}\mathbb{E}[X_1].
\end{equation*}
With the help of (\ref{mv1}), we get
\begin{align*}
\text{Var}[C^{\beta}(t)] 
&= \mathbb{E}[\mathcal{N}^{\beta}(t)] \text{Var}[X_1] + \text{Var}[\mathcal{N}^{\beta}(t)] (\mathbb{E}[X_1])^2\\
&= St^{\beta} \text{Var}[X_1] + (Rt^{2 \beta}+ T t^{\beta}) (\mathbb{E}[X_1])^2\\
&= St^{\beta}\mathbb{E}[X_1^2] - St^{\beta}(\mathbb{E}[X_1])^2 + (Rt^{2 \beta}+ T t^{\beta}) (\mathbb{E}[X_1])^2\\
&= St^{\beta}\mathbb{E}[X_1^2] + (Rt^{2 \beta}+ (T-S ) t^{\beta}) (\mathbb{E}[X_1])^2.
\end{align*}

\begin{remark}
	It is observed that $\text{Var}[C^{\beta}(t)] - \mathbb{E}[C^{\beta}(t)] > 0$ for all $t \geq 0,$ since $E(X_1^2) - E(X_1) \geq 0$ as Pr$\{X_1 \geq 1\} = 1$.  Hence, the CGFCP is over-dispersed. 
\end{remark}
\noindent In the next proposition, we discuss a particular case of the CGFCP.
\begin{proposition}
	Let $X_j,\; j=1,2 \ldots$ be sequence of i.i.d. random variables with pmf Pr$\{X_1= j\} =  \beta_j, \;  j=0,1,2, \ldots.$ Then, the probability generating function (pgf) $\mathcal{G}^{\beta}(u, t)$ of the process $Q(t) = \sum_{j=1}^{\mathcal{N}^{\beta}(t)} X_j $ is given by
	\begin{equation*}
	\mathcal{G}^{\beta}(u, t) = E_{\beta, 1}^1 \left(t^{\beta} \sum_{i=1}^{k}\lambda_i \sum_{j=0}^\infty \beta_j^{*i}(u^j - 1) \right),
	\end{equation*}
	where $\beta_j^{*i} = \sum_{\substack{
			l_1+ l_2+ \cdots+ l_i =j \\ \\ l_i \in \mathbb{N}_0} } \beta_{l_1}  \beta_{l_2}\cdots  \beta_{l_i}$ for all $i=1,2, \dots k.$
\end{proposition}
\begin{proof}
	The pgf of the GFCP is given by \cite[formula (14)]{Kataria2022}
	\begin{equation}\label{pgf1}
	\mathbb{E}[u^{\mathcal{N}^{\beta}(t)}] = E_{\beta, 1}^1 \left(\sum_{i=1}^{k}\lambda_i (u^i -1)t^{\beta}\right),\;\; |u| \leq 1.
	\end{equation}
	We define $G_{X_1}^{*k}(u) = \mathbb{E}[u^{{\sum_{j=1}^k}X_j}]$. Then, from (\ref{pgf1}), we have
	\begin{align*}
	\mathcal{G}^{\beta}(u, t) 
	&= \mathbb{E}[u^{Q(t)}] = E_{\beta, 1}^1 \left(\sum_{i=1}^{k}\lambda_i \left((G_{X_1}(u))^i -1\right )t^{\beta}\right)\\
	&= E_{\beta, 1}^1 \left(t^{\beta}\sum_{i=1}^{k}\lambda_i \prod_{j=1}^{k}G_{X_j}(u) - \Delta t^{\beta}\right)\\
	&= E_{\beta, 1}^1 \left(t^{\beta} \sum_{i=1}^{k}\lambda_i G_{X_1+X_2+\cdots+X_k}(u) - \Delta t^{\beta}\right)\\
	&= E_{\beta, 1}^1 \left(t^{\beta} \sum_{i=1}^{k}\lambda_i \sum_{j=0}^\infty \text{Pr}\{X_1+X_2+\cdots+ X_k = j\}u^j - \Delta t^{\beta}\right)\\
	&= E_{\beta, 1}^1 \left(t^{\beta} \sum_{i=1}^{k}\lambda_i \sum_{j=0}^\infty \beta_j^{*i}u^j - \Delta t^{\beta}\right)
	= E_{\beta, 1}^1 \left(t^{\beta} \sum_{i=1}^{k}\lambda_i \sum_{j=0}^\infty \beta_j^{*i}(u^j - 1) \right).
	\end{align*}
	Hence, the proposition follows.
\end{proof}
		
\section{Generalized Fractional Risk Process}
 In this section, we present a generalized fractional risk model such that the claim count process is observed to be a GFCP. We call it generalized fractional risk process (GFRP). This model generalizes the fractional risk process based on the fractional Poisson process and provides a more general and a flexible framework. For $\beta \in (0,1)$, we consider a risk process $R^{\beta}(t)$ defined as
\begin{equation}\label{rm1}
R^{\beta}(t) = \nu + \eta (1+ \rho) \sum_{i=1}^{k}i \lambda_i \mathcal{W}_\beta(t) - \sum_{j=1}^{\mathcal{N}^{\beta}(t)} X_j,
\end{equation}			
where $\nu >0$ is the initial capital, $X_j, \; j \geq 1$ are i.i.d. random claim size independent of $\mathcal{N}^{\beta}(t)$ with $\mathbb{E}[X_j]= \eta  >0$. The parameter $\rho \geq 0$ is the safety loading parameter which is essential for the GFRP for satisfying the net profit condition.

In the next result, we show that the generalized fractional risk model in (\ref{rm1}) satisfies the martingale conditions depending on the safety loading parameter $\rho$. First, we proceed with the following Lemma from \cite{Kataria2021}.
\begin{lemma}\label{lm1}
	Let $\mathcal{N}(t)$ be a generalized counting process with positive rates $\lambda_1, \lambda_2,\ldots, \lambda_k.$ Then, the process $\left\{\mathcal{N}(t) -  \sum_{i=1}^{k}i \lambda_i t \right\}_{t \geq 0}$ is a martingale with respect to a natural filtration $\mathcal{F}_t =\sigma (\mathcal{N}(s), \; s \leq t). $
\end{lemma}

\begin{theorem}\label{thm1}
	The generalized fractional risk process in (\ref{rm1}) is a martingale (submartingale, and supermartingale) for $\rho =0 ~(\rho >0, \text{ and } \rho < 0),$ respectively with respect to the filtration $\mathcal{F}_t = \sigma\left(\mathcal{N}^{\beta}(s), \; s \leq t\right) \vee \sigma\left(\mathcal{W}_\beta(s), \; s \geq 0\right).$
\end{theorem}

\begin{proof}
	For $0 < s \leq t$, we have
	\begin{align*}
\mathbb{E}[R^{\beta}(t)-R^{\beta}(s) \; | \; \mathcal{F}_s]	&= \mathbb{E}\left[\eta (1+ \rho) \sum_{i=1}^{k}i \lambda_i \mathcal{W}_\beta(t) - \sum_{j=1}^{\mathcal{N}^{\beta}(t)} X_j - \left(\eta (1+ \rho) \sum_{i=1}^{k}i \lambda_i \mathcal{W}_\beta(s) - \sum_{j=1}^{\mathcal{N}^{\beta}(s)} X_j\right)\; \vline \; \mathcal{F}_s\right]\\
	&=  \mathbb{E}\left[\eta (1+ \rho) \sum_{i=1}^{k}i \lambda_i \left(\mathcal{W}_\beta(t) - \mathcal{W}_\beta(s)\right) - \sum_{j=\mathcal{N}^{\beta}(s)+1}^{\mathcal{N}^{\beta}(t)} X_j \; \vline \; \mathcal{F}_s\right]\\
	&=  \mathbb{E}\left[\eta (1+ \rho) \sum_{i=1}^{k}i \lambda_i \left(\mathcal{W}_\beta(t) - \mathcal{W}_\beta(s)\right)\; \vline \; \mathcal{F}_s\right]- \mathbb{E}\left[\mathbb{E}\left[ \sum_{j=\mathcal{N}^{\beta}(s)+1}^{\mathcal{N}^{\beta}(t)} X_j\; \vline \; \mathcal{F}_t\right]\; \vline \; \mathcal{F}_s\right]\\
	&=  \mathbb{E}\left[\eta (1+ \rho) \sum_{i=1}^{k}i \lambda_i \left(\mathcal{W}_\beta(t) - \mathcal{W}_\beta(s)\right)\; \vline \; \mathcal{F}_s\right]- \mathbb{E}\left[\eta \left(\mathcal{N}^{\beta}(t) - \mathcal{N}^{\beta}(s)\right)\; \vline \; \mathcal{F}_s\right]\\
	&=  \mathbb{E}\left[\eta (1+ \rho) \sum_{i=1}^{k}i \lambda_i \left(\mathcal{W}_\beta(t) - \mathcal{W}_\beta(s)\right) - \eta \left(\mathcal{N}^{\beta}(t) - \mathcal{N}^{\beta}(s)\right)\; \vline \; \mathcal{F}_s\right]\\
	&= -\eta \mathbb{E}\left[\left(\mathcal{N}^{\beta}(t) -  \sum_{i=1}^{k}i \lambda_i \mathcal{W}_\beta(t) \right)- \left(\mathcal{N}^{\beta}(s) -  \sum_{i=1}^{k}i \lambda_i \mathcal{W}_\beta(s) \right) \; \vline \; \mathcal{F}_s\right] 
	+ \eta \rho  \sum_{i=1}^{k}i \lambda_i \mathbb{E}\left[\mathcal{W}_\beta(t) - \mathcal{W}_\beta(s) \; \vline \; \mathcal{F}_s\right]\\
	&=  \eta \rho  \sum_{i=1}^{k}i \lambda_i \mathbb{E}\left[\mathcal{W}_\beta(t) - \mathcal{W}_\beta(s) \; \vline \; \mathcal{F}_s\right].
	\end{align*}
	It is essential to observe that the claim count process and the payment process $W_\beta(t)$  in the model (\ref{rm1}) follow the same time scale. The net and equivalence principle (see \cite{Mikosch2009}) is consistent if we time-change the classical risk model by the inverse $\beta$-stable subordinator considering Lemma \ref{lm1} as $\left\{\sum_{j=1}^{\mathcal{N}^{\beta}(t)} X_j - \eta  \sum_{i=1}^{k}i \lambda_i \mathcal{W}_\beta(t) \right\}_{t \geq 0}$ is a martingale. Hence, the last step in the proof  follows. 
\end{proof}

\subsection{Variants of the Risk Process in (\ref{rm1})} Recently, risk models based on the non-homogeneous FPP and the non-homogeneous mixed fractional Poisson process have been considered in the literature (see \cite{KatariaRisk} and \cite{Kumar2020}). We here consider a variant of the GFRP defined as
\begin{equation}\label{nhrp1}
\hat{R}^{\beta}(t) = \nu + \eta (1+ \rho) \Lambda(\mathcal{W}_\beta(t)) - \sum_{j=1}^{\mathcal{N}_{\Lambda}^{\beta}(t)} X_j,
\end{equation}
where $\Lambda$ is a right continuous non-decreasing function with $\Lambda(0) =0, \Lambda(\infty) =\infty,$ and $\Lambda(t)- \Lambda(t-) \leq 1$.

\begin{remark}
	When $ \Lambda(t) = t \sum_{i=1}^{k}i \lambda_i$, the risk process $	\hat{R}^{\beta}(t)$ leads to GFRP.
\end{remark}

\begin{lemma}\label{lm2}
	Let $\mathcal{N}_\lambda(t)$ be a non-homogeneous version of the generalized counting process with positive rates $\lambda_1, \lambda_2,\ldots, \lambda_k.$ Then, the process $\left\{\mathcal{N}_\lambda(t) -  \sum_{i=1}^{k}i \lambda_i \Lambda(t) \right\}_{t \geq 0}$ is a martingale with respect to a natural filtration $\mathcal{F}_t =\sigma (\mathcal{N}_\lambda(s), \; s \leq t). $
\end{lemma}
\begin{proof}
	Since the non-homogeneous version of the GCP has the independent increment property, the lemma follows. 
\end{proof}
\begin{theorem}
	The risk process in (\ref{nhrp1}) is a martingale (submartingale, and supermartingale) for $\rho =0~(\rho >0~ \text{and} ~\rho < 0),$ respectively with respect to the filtration $\mathcal{F}_t = \sigma\left(\mathcal{N}_\Lambda^{\beta}(s), \; s \leq t\right) \vee \sigma\left(\mathcal{W}_\beta(s), \; s \geq 0\right).$
\end{theorem}
\begin{proof}
	The proof of the theorem follows on the same line as proof of Theorem \ref{thm1} with the help of Lemma \ref{lm2}.
\end{proof}
Similarly, we may provide several other variants of the risk process by considering the expectation of the inverse $\beta$-stable subordinator and $\mathbb{E}[\Lambda(\mathcal{W}_\beta(t))].$
\subsection{Dependence Structure of the GFRP in (\ref{rm1})} \hfill \\
In this section, we study dependence structure of the GFRP.
\begin{theorem}\label{Cov1}
	For $0 <  s \leq t$, the covariance of the GFRP is given by
	\begin{equation}\label{cov111}
	\text{Cov}[R^{\beta}(t), R^{\beta}(s)] = \eta^2 \left(\sum_{j=1}^{k}j \lambda_j\right)^2 \rho^2 \text{Cov}[\mathcal{W}_\beta(s), \mathcal{W}_\beta(t)]+ Ss^{\beta} \text{Var}[X_i] + Ts^{\beta}\eta^2.
	\end{equation}
\end{theorem}
\begin{proof} Let $I\{\cdot\}$ denotes the indicator function.
For a counting phenomena, we have $\mathcal{N}^{\beta}(s) \leq \mathcal{N}^{\beta}(t)$ for $s \leq t$.\\ Consider
	\begin{align*}
	\text{Cov}\left[\sum_{j=1}^{\mathcal{N}^{\beta}(s)} X_j, \sum_{i=1}^{\mathcal{N}^{\beta}(t)} X_i\right]&= \mathbb{E}\left[\sum_{j=1}^{\infty} \sum_{i=1}^{\infty}X_j X_i \;I\{\mathcal{N}^{\beta}(t) \geq i, \mathcal{N}^{\beta}(s) \geq j\}\right] - (\mathbb{E}[X_i])^2 \mathbb{E}[\mathcal{N}^{\beta}(t)] \mathbb{E}[\mathcal{N}^{\beta}(s)]\\
	&= \mathbb{E}\left[\sum_{j=1}^{\infty} X_j^2 \; I\{\mathcal{N}^{\beta}(s) \geq j\}\right] +  \mathbb{E}\left[\mathop{\sum^{}\sum^{}}_{\ i \neq j\ }X_j X_i \;I\{\mathcal{N}^{\beta}(t) \geq i, \mathcal{N}^{\beta}(s) \geq j\}\right]\\
	& \;\;\;  - (\mathbb{E}[X_i])^2 \mathbb{E}[\mathcal{N}^{\beta}(t)] \mathbb{E}[\mathcal{N}^{\beta}(s)]\\
	&= \mathbb{E}[X_j^2] \sum_{j=1}^{\infty} \text{Pr}\{\mathcal{N}^{\beta}(s) \geq j\} + (\mathbb{E}[X_i])^2 \sum_{j=1}^{\infty} \sum_{i=1}^{\infty} \text{Pr}\{\mathcal{N}^{\beta}(t) \geq i, \mathcal{N}^{\beta}(s) \geq j\} \\
	& \;\;\; - (\mathbb{E}[X_i])^2\sum_{j=1}^{\infty} \text{Pr}\{\mathcal{N}^{\beta}(s) \geq j\}
	- (\mathbb{E}[X_i])^2 \mathbb{E}[\mathcal{N}^{\beta}(t)] \mathbb{E}[\mathcal{N}^{\beta}(s)]\\
	&= \mathbb{E}[X_i^2] \mathbb{E}[\mathcal{N}^{\beta}(s) ]+ (\mathbb{E}[X_i])^2 \mathbb{E}[\mathcal{N}^{\beta}(s) \mathcal{N}^{\beta}(t)]\\
	&\;\;\; - (\mathbb{E}[X_i])^2 \mathbb{E}[\mathcal{N}^{\beta}(s) ] - (\mathbb{E}[X_i])^2 \mathbb{E}[\mathcal{N}^{\beta}(t)] \mathbb{E}[\mathcal{N}^{\beta}(s)].
	\end{align*}
	With the help of (\ref{mv1}) and (\ref{cov11}), we get
	\begin{align}\label{DS1}
	\text{Cov}\left[\sum_{j=1}^{\mathcal{N}^{\beta}(s)} X_j, \sum_{i=1}^{\mathcal{N}^{\beta}(t)} X_i\right] 
	&= \text{Var}[X_i]\mathbb{E}[\mathcal{N}^{\beta}(s) ] + (\mathbb{E}[X_i])^2 \text{Cov}[\mathcal{N}^{\beta}(s), \mathcal{N}^{\beta}(t)]\nonumber\\
	&= Ss^{\beta} \text{Var}[X_i] + (\mathbb{E}[X_i])^2 \left(Ts^{\beta} + \left(\sum_{j=1}^{k}j \lambda_j\right)^2 \text{Cov}[\mathcal{W}_\beta(s), \mathcal{W}_\beta(t)]\right).
	\end{align}
	Also, we have
	\begin{align}\label{DS2}
	\text{Cov}\left[\mathcal{W}_\beta(s), \sum_{i=1}^{\mathcal{N}^{\beta}(t) }X_i\right] 
	&= \mathbb{E}\left[\mathcal{W}_\beta(s) \sum_{i=1}^{\mathcal{N}^{\beta}(t) }X_i\right] - \mathbb{E}[\mathcal{W}_\beta(s)] \mathbb{E}\left[\sum_{i=1}^{\mathcal{N}^{\beta}(t) }X_i\right]\nonumber\\
	&= \mathbb{E}\left[\mathcal{W}_\beta(s)\mathbb{E}\left[ \sum_{i=1}^{\mathcal{N}^{\beta}(t) }X_i\; \vline \; \mathcal{W}_\beta(s), \mathcal{N}^{\beta}(t) \right]\right]	- \mathbb{E}[X_i]\mathbb{E}[\mathcal{N}^{\beta}(t)]\mathbb{E}[\mathcal{W}_\beta(s)]\nonumber\\
	&= \mathbb{E}[X_i] \mathbb{E}[\mathcal{W}_\beta(s)W_\beta(t)]\sum_{j=1}^{k}j \lambda_j - \mathbb{E}[X_i]\mathbb{E}[\mathcal{N}^{\beta}(t)]\mathbb{E}[\mathcal{W}_\beta(s)]\nonumber\\
	&= \mathbb{E}[X_i] \text{Cov}[\mathcal{W}_\beta(s), \mathcal{W}_\beta(t)] \sum_{j=1}^{k}j \lambda_j.
	\end{align}
	Similarly, we have
	\begin{align}\label{DS3}
	\text{Cov}\left[\sum_{i=1}^{\mathcal{N}^{\beta}(s) }X_i, \mathcal{W}_\beta(t)\right] = \mathbb{E}[X_i] \text{Cov}[\mathcal{W}_\beta(s), \mathcal{W}_\beta(t)] \sum_{j=1}^{k}j \lambda_j.
	\end{align}
	When $\mathbb{E}[X_1] = \eta,$ with the help of (\ref{DS1}), (\ref{DS2}) and (\ref{DS3}) the covariance of the GFRP is derived as\\
	\begin{align*}
\text{Cov}[R^{\beta}(t), R^{\beta}(s)]	&= \eta^2 \left(\sum_{j=1}^{k}j \lambda_j\right)^2 (1+\rho)^2 \text{Cov}[\mathcal{W}_\beta(s), \mathcal{W}_\beta(t)] - \eta \sum_{j=1}^{k}j \lambda_j 	\text{Cov}\left[\mathcal{W}_\beta(s), \sum_{i=1}^{\mathcal{N}^{\beta}(t) }X_i\right] \\
	&\;\;\; + \eta \sum_{j=1}^{k}j \lambda_j \text{Cov}\left[\sum_{i=1}^{\mathcal{N}^{\beta}(s) }X_i, \mathcal{W}_\beta(t),\right] + \text{Cov}\left[\sum_{j=1}^{\mathcal{N}^{\beta}(s)} X_j, \sum_{i=1}^{\mathcal{N}^{\beta}(t)} X_i\right]\\
	&= \eta^2 \left(\sum_{j=1}^{k}j \lambda_j\right)^2 \rho^2 \text{Cov}[\mathcal{W}_\beta(s), \mathcal{W}_\beta(t)]+ Ss^{\beta} \text{Var}[X_i] + Ts^{\beta}\eta^2.
	\end{align*}
\end{proof}

\begin{remark} When $s=t$, (\ref{cov111}) yields the covariance of GFRP of the following form
	\begin{align*}
	\text{Var}[R^{\beta}(t)] 
	&= \eta^2 \left(\sum_{j=1}^{k}j \lambda_j\right)^2 \rho^2 \text{Var}[ \mathcal{W}_\beta(t)]+ St^{\beta} \text{Var}[X_i] + Tt^{\beta}\eta^2.
	\end{align*}
\end{remark}

\begin{proposition}
	For $\beta \in (0,1)$, the GFRP $R^{\beta}(t)$ exhibits the LRD property.
\end{proposition}
\begin{proof}
	By observing the asymptotic behaviour of the covariance and the variance of the inverse $\beta$-stable subordinator $\mathcal{W}_\beta(t)$ studied in \cite[Section 2.2(v)]{Kumar2020}, we obtain the correlation function of the GFRP as
	\begin{align*}
	\text{Corr}[R^{\beta}(t), R^{\beta}(s)] 
	&= \frac{\text{Cov}[R^{\beta}(t), R^{\beta}(s)]}{\sqrt{\text{Var}[R^{\beta}(s)] }\sqrt{\text{Var}[R^{\beta}(t)] }}\\
	& \sim t^{-\beta}\frac{ \eta^2 \left(\sum_{j=1}^{k}j \lambda_j\right)^2 \rho^2 B(1+\beta, \beta)s^{2\beta} + Ss^{\beta} \text{Var}[X_i] + Ts^{\beta}\eta^2}{\eta \rho \sqrt{d(\beta)} \sqrt{\text{Var}[R^{\beta}(s)] } \sum_{j=1}^{k}j \lambda_j },
	\end{align*}
	where $d(\beta) = \frac{2}{\Gamma(1+2 \beta)} - \frac{1}{\Gamma^2(1+\beta)}$ and $B(1+\beta, \beta)$ is the beta function. Since $\beta \in (0,1)$. By LRD definition, the GFRP has the LRD property.
\end{proof}

Next, we introduce the generalized fractional noise risk process and discuss its short-range dependence (SRD) property.
\subsection{Generalized Fractional Noise Risk Process} \hfill \\
For $\varepsilon > 0$, the generalized fractional noise risk process (GFNRP) is the increment $H_\varepsilon^{\beta}(t)$ defined as
\begin{equation*}
H_\varepsilon^{\beta}(t) = R^{\beta}(t+ \varepsilon)- R^{\beta}(t).
\end{equation*}
\begin{theorem}
	For $\beta \in (0,1)$, the GFNRP $H_\varepsilon^{\beta}(t)$ exhibits the SRD property.
\end{theorem}
\begin{proof}
	For $ 0 < s +\varepsilon \leq t$, from (\ref{cov111}), we have
	\begin{align*}
	\text{Cov}[H_\varepsilon^{\beta}(s), H_\varepsilon^{\beta}(t)] 
	&= \text{Cov}[R^{\beta}(t+\varepsilon), R^{\beta}(s+\varepsilon)]+ \text{Cov}[R^{\beta}(t), R^{\beta}(s)]\\
	&\;\;\;- \text{Cov}[R^{\beta}(t+\varepsilon), R^{\beta}(s)]- \text{Cov}[R^{\beta}(t), R^{\beta}(s+ \varepsilon)]\\
	&= \eta^2 \left(\sum_{j=1}^{k}j \lambda_j\right)^2 \rho^2 \{\text{Cov}[\mathcal{W}_\beta(s), \mathcal{W}_\beta(t)] - \text{Cov}[\mathcal{W}_\beta(s+\varepsilon), \mathcal{W}_\beta(t)]\\
	&\;\;\; - \text{Cov}[\mathcal{W}_\beta(s), \mathcal{W}_\beta(t+\varepsilon)]+\text{Cov}[\mathcal{W}_\beta(s+\varepsilon), \mathcal{W}_\beta(t+\varepsilon)]\}.
	\end{align*}

	Therefore, with the help of (\ref{asy11}), we get
	\begin{align*}
	\text{Cov}[H_\varepsilon^{\beta}(s), H_\varepsilon^{\beta}(t)] 
	&\sim \eta^2 \left(\sum_{j=1}^{k}j \lambda_j\right)^2 \rho^2  \frac{1}{\Gamma^2(\beta +1)} \left(\frac{\beta^2}{\beta+1}\right)\\
	&\;\;\; \times \left[\frac{(s+\varepsilon)^{\beta+1}}{t^{1-\beta}} + \frac{s^{\beta+1}}{(t+\varepsilon)^{1-\beta}}- \frac{s^{\beta+1}}{t^{1-\beta}}- \frac{(s+\varepsilon)^{\beta+1}}{(t+\varepsilon)^{1-\beta}}\right]\\
	&\sim \eta^2 \left(\sum_{j=1}^{k}j \lambda_j\right)^2 \rho^2  \frac{1}{\Gamma^2(\beta +1)} \left(\frac{\beta^2}{\beta+1}\right) \varepsilon (1-\beta) \left((s+\varepsilon)^{\beta+1}-s^{\beta+1} \right)t^{\beta-2}.
	\end{align*}
	Also, we have
	\begin{align*}
	\text{Var}[H_\varepsilon^{\beta}(t)] 
	&= 	\text{Var}[R^{\beta}(t+ \varepsilon)] + 	\text{Var}[R^{\beta}(t)] - 2 	\text{Cov}[R^{\beta}(t), R^{\beta}(t+\varepsilon)] \\
	&= \eta^2 \left(\sum_{j=1}^{k}j \lambda_j\right)^2 \rho^2 \text{Var}[ \mathcal{W}_\beta(t+\varepsilon)]+ 2S(t+\varepsilon)^{\beta} \text{Var}[X_i] + 2T(t+\varepsilon)^{\beta}\eta^2\\
	&\;\;\;+  \eta^2 \left(\sum_{j=1}^{k}j \lambda_j\right)^2 \rho^2 \text{Var}[ \mathcal{W}_\beta(t)]+ St^{\beta} \text{Var}[X_i] + Tt^{\beta}\eta^2\\
	&\;\;\;- \eta^2 \left(\sum_{j=1}^{k}j \lambda_j\right)^2 \rho^2 \text{Cov}[\mathcal{W}_\beta(t+\varepsilon), \mathcal{W}_\beta(t)].
	\end{align*}
	Since $B(\beta, \beta+1; t/(t+\varepsilon)) \sim B(\beta, \beta+1)$ for large $t$ and let $d(\beta) = \left(\frac{2}{\Gamma(2 \beta+1)}- \frac{1}{\Gamma^2(1+\beta)}\right)$. Then, with the help of (\ref{asy11}) and formula (11) in \cite{Leonenko2014}, we get
	\begin{align*}
	\text{Var}[H_\varepsilon^{\beta}(t)] 
	& \sim \eta^2 \left(\sum_{j=1}^{k}j \lambda_j\right)^2 \rho^2 d(\beta) (t^{2 \beta} + (t+\varepsilon)^{2 \beta})- S(t^{\beta}-(t+\varepsilon)^{\beta} )\text{Var}[X_i] + T(t^{\beta}- (t+\varepsilon)^{\beta} )\eta^2\\
	&\;\;\;- \left(\sum_{j=1}^{k}j \lambda_j\right)^2 \frac{ 2\eta^2\rho^2 }{\Gamma^2(\beta+1)} \left(\beta t^{2 \beta} B(\beta, \beta+1) - \beta (t + \varepsilon)^{2 \beta}B\left(\beta, \beta+1; \frac{t}{t+\varepsilon}\right) +(t(t+\varepsilon))^{\beta} \right)\\
	&\sim \eta^2 \varepsilon\left(\sum_{j=1}^{k}j \lambda_j\right)^2 \rho^2 \left(\frac{2}{\Gamma(2 \beta+1)}- \frac{1}{\Gamma^2(1+\beta)}\right) t^{2 \beta} + (T \eta^2 -S \text{Var}[X_i] )\varepsilon t^{\beta-1}\\
	&\;\;\;- \left(\sum_{j=1}^{k}j \lambda_j\right)^2 \frac{ 2\eta^2\rho^2 }{\Gamma^2(\beta+1)} \left(\beta t^{2 \beta} B(\beta, \beta+1) - \beta (t + \varepsilon)^{2 \beta}B\left(\beta, \beta+1\right) +(t(t+\varepsilon))^{\beta} \right)\\
	&\sim (T \eta^2 -S \text{Var}[X_i] ) t^{\beta-1}.
	\end{align*}
	Therefore, the correlation function of the GFNRP takes the following form
	\begin{align*}
	\text{Corr}[H_\varepsilon^{\beta}(t), H_\varepsilon^{\beta}(s)] 
	&= \frac{\text{Cov}[H_\varepsilon^{\beta}(t), H_\varepsilon^{\beta}(s)]}{\sqrt{\text{Var}[H_\varepsilon^{\beta}(t)] }\sqrt{\text{Var}[H_\varepsilon^{\beta}(s)] }} \sim B(s) t^{-(3-\beta)/2}, \;\; \text{as } t \rightarrow \infty,
	\end{align*}
	where $$B(s) = \frac{\eta^2 \left(\sum_{j=1}^{k}j \lambda_j\right)^2 \rho^2  \frac{1}{\Gamma^2(\beta +1)} \left(\frac{\beta^2}{\beta+1}\right) \varepsilon (1-\beta) \left((s+\varepsilon)^{\beta+1}-s^{\beta+1} \right)}{\sqrt{ (T \eta^2 -S \text{Var}[X_i] )\;\text{Var}[H_\varepsilon^{\beta}(s)] } }.$$
	It is clear that $\beta \in (0,1) \Rightarrow (3-\beta)/2 \in (1, 3/2)$ that proves the SRD property.
\end{proof}

\section{An Alternative to Generalized Fractional Risk process}
\noindent In this section, we give an alternative to the GFRP. Let $\mathcal{N}^{\beta}(t)$ be the GFCP with stability index $\beta \in (0,1)$. Then, we define the risk process
\begin{equation}\label{rpa1}
\underline{R}^{\beta}(t) = \nu + \eta t - \sum_{j=1}^{\mathcal{N}^{\beta}(t)}Y_j, \;\; t \geq 0,
\end{equation}
where $Y_j, \;j=1,2,\ldots$ are i.i.d. random variables with mean $\mu>0$ and finite variance and are independent of $\mathcal{N}^{\beta}(t) = \mathcal{N}^{\beta}(S_\beta(t))$, and  $\eta > 0$ is the constant premium rate. We abbreviate it as AGFRP. The dependence structure of the AGFRP is similar to the GFRP. 

Let $\Psi(u,t)$ be the probability of ruin in finite time defined as
\begin{equation*}
\Psi(u,t) = \text{Pr} \{\underline{R}^{\beta}(s) <0 \text{ for some } s \leq t\}.
\end{equation*}
Next, we discuss the ruin probabilities of the risk process with light-tailed and heavy-tailed claim sizes under the assumptions of the model in (\ref{rpa1}).

\subsection{Ruin Probability with Light-Tailed Claim Sizes} \hfill \\
With the aim of providing asymptotic results for the ruin probability $\Psi(u,t)$, we first recall the CGFCP defined as
\begin{equation*}
C^{\beta}(t) = \sum_{j=1}^{\mathcal{N}^{\beta}(t)} Y_j,\;\; t \geq 0,
\end{equation*}
where $Y_j, \; j=1,2,\ldots$ are i.i.d. random variables with common moment generating function (MGF) $H_Y$ and are independent of $\mathcal{N}^{\beta}(t)$. Hence, the mgf of the CGFCP is derived as
\begin{align*}
\mathbb{E}[e^{\xi 	C^{\beta}(t) }] = 	\mathbb{E}\left[\mathbb{E}\left[e^{\xi 	\sum_{j=1}^{\mathcal{N}^{\beta}(t)} Y_j}\; \vline \; \mathcal{N}^{\beta}(t)\right]\right] = \mathbb{E}\left[\mathbb{E}\left[\left(H_Y(\xi)\right)^{\mathcal{N}^{\beta}(t)}\; \vline \; \mathcal{N}^{\beta}(t)\right]\right].
\end{align*}
With the help of formula (14) in \cite{Kataria2021}, we get
\begin{equation}\label{mgf1}
\mathbb{E}[e^{\xi 	C^{\beta}(t) }] = E_{\beta, 1}^1\left(\sum_{i=1}^k \lambda_i \left(\left(H_Y(\xi)\right)^i -1\right)t^{\beta}\right).
\end{equation}
Next, in the following proposition, we present asymptotic structure of the ruin probability in finite time.
\begin{proposition}
		Under the assumption of the risk model in (\ref{rpa1}) when the size of succesive claims $Y_j$'s are light-tailed, we have
		\begin{equation}\label{ine1}
		\frac{1}{2}\Psi(u,t) \leq e^{-u \xi_0(t, \beta, \lambda_1, \lambda_2, \ldots, \lambda_k)},
		\end{equation}
		where $\xi_0(t, \beta, \lambda_1, \lambda_2, \ldots, \lambda_k) $ is a unique real number with $t >0$.
	\end{proposition}
	\begin{proof}
		The proof is motivated by proof of Proposition 6 in \cite{Biard2014}. Using Chebyshev inequality and with the help of (\ref{mgf1}), we get
		\begin{equation}\label{asy1}
		\text{Pr} \left\{\sum_{j=1}^{\mathcal{N}^{\beta}(t)} Y_j > u\right \} = \text{Pr} \left\{e^{\xi 	C^{\beta}(t) } > e^{\xi u}\right\} \leq e^{-\xi u} E_{\beta, 1}^1\left(\sum_{i=1}^k \lambda_i \left(\left(H_Y(\xi)\right)^i -1\right)t^{\beta}\right).
		\end{equation}
		Since for $|y| < 1$, we have (see \cite[Formula (22)]{Biard2014})
		\begin{equation}\label{asymit1}
		E_{\beta, 1}^1(y) \leq \frac{1}{\Gamma(y_0)}\sum_{k=0}^\infty y^k \leq \frac{1}{(1-y)\Gamma(y_0)},
		\end{equation}
		where $\Gamma(y_0) \simeq 0.8856$ to achieve minimum of the Mittag-Leffler function (for more details see Proposition 6 in \cite{Biard2014}).\\
		We also define a unique positive real number $\xi_0$ such that	$$\sum_{i=1}^k \lambda_i (H_Y(\xi_0))^i = 1+ \frac{2 \Gamma(y_0)-1}{2t^\beta \Gamma(y_0) }.$$
		Therefore, from (\ref{asy1}) and (\ref{asymit1}), we get
		\begin{align*}
		\text{Pr} \left\{\sum_{j=1}^{\mathcal{N}^{\beta}(t)} Y_j > u\right \} \leq \frac{e^{-\xi_0 u}}{\Gamma(y_0)\left(1- \sum_{i=1}^k \lambda_i \left(\left(H_Y(\xi_0)\right)^i -1\right)t^{\beta}\right)} \leq e^{-u \xi_0}.
		\end{align*}
		Hence, the inequality (\ref{ine1})  is established.
\end{proof}
\subsection{Ruin Probability with Heavy-Tailed Claim Sizes} \hfill \\
Let $Y_j, \; j=1,2, \ldots$ be subexponentially distributed i.i.d. random varaibles with distribution function $G_Y(t) =$Pr$\{Y_j \leq t\}$ satisfying $$\lim_{t \rightarrow \infty} \left(1-G_Y^{*2}(t)\right)/ \left(1-G_Y(t)\right) =2,$$
where $G_Y^{*2}$ denotes the two-fold convolution of $G_Y$.
\begin{proposition}
	Under the assumption of the risk model in (\ref{rpa1}) if the sizes of successive claims $Y_i$'s are subexponentially distributed with survival function $\bar{G}_Y = 1-G_Y$, then the ruin probability satisfies
	\begin{equation*}
	\Psi(u,t) \sim  \frac{t^\beta}{\Gamma(\beta +1)}\sum_{i=1}^k i \lambda_i\bar{G}_Y(u), \text{ as } u \rightarrow \infty.
	\end{equation*}
\end{proposition}
\begin{proof}
	For the proof one can refer to the  Proposition 4 of \cite{Biard2014}.
\end{proof}

\bibliographystyle{apalike}
\bibliography{rmb}
\end{document}